\documentclass{article}
\usepackage{amsmath,amsthm,bm,colortbl,graphicx,mathrsfs,afterpage,amssymb,url}
\theoremstyle{plain}
\newtheorem{theorem}{Theorem}
\newtheorem{lemma}[theorem]{Lemma}

\theoremstyle{remark}
\def\P{{\rm P}} 
\def\E{{\rm E}} 
\def\Var{{\rm Var}} 
\def\Cov{{\rm Cov}} 
\def\T{\textsf{T}} 
\def\eps{\varepsilon} 
\allowdisplaybreaks
\makeatletter 
\def\@biblabel#1{\hspace*{-\labelsep}}
\makeatother
\title{Parrondo's paradox via redistribution of wealth}
\author{S. N. Ethier\thanks{This work was partially supported by a grant from the Simons Foundation (209632 to S.N.E.).  It was also supported by the Korean Federation of Science and Technology Societies Grant funded by the Korean Government (MEST, Basic Research Promotion Fund).}\\
\begin{small}University of Utah\end{small}\\
\begin{small}Department of Mathematics\end{small}\\
\begin{small}155 S. 1400 E., JWB 233\end{small}\\
\begin{small}Salt Lake City, UT 84112, USA\end{small}\\
\begin{small}e-mail: \url{ethier@math.utah.edu}\end{small}
\and
Jiyeon Lee\thanks{This work was supported by the Basic Science Research Program through the National Research Foundation of Korea (NRF) funded by the Ministry of Education, Science and Technology (2011-0005982).}\\
\begin{small}Yeungnam University\end{small}\\
\begin{small}Department of Statistics\end{small}\\
\begin{small}214-1 Daedong, Kyeongsan\end{small}\\
\begin{small}Kyeongbuk 712-749, South Korea\end{small}\\
\begin{small}e-mail: \url{leejy@yu.ac.kr}\end{small}}
\date{}

\date{Sept.\ 20, 2011}

\begin{document}
\maketitle

\begin{abstract}
Toral (2002) considered an ensemble of $N\ge2$ players.  In game $B$ a player is randomly selected to play Parrondo's original capital-dependent game.  In game $A'$ two players are randomly selected without replacement, and the first transfers one unit of capital to the second.  Game $A'$ is fair (with respect to total capital), game $B$ is losing (or fair), and the random mixture $\gamma A'+(1-\gamma)B$ is winning, as was demonstrated by Toral for $\gamma=1/2$ using computer simulation.  We prove this, establishing a strong law of large numbers and a central limit theorem for the sequence of profits of the ensemble of players for each $\gamma\in(0,1)$.  We do the same for the nonrandom pattern of games $(A')^r B^s$ for all integers $r,s\ge1$.  An unexpected relationship between the random-mixture case and the nonrandom-pattern case occurs in the limit as $N\to\infty$.\medskip\par
\end{abstract}

\newpage
\section{Introduction}

The original Parrondo (1996) games can be described in terms of probabilities $p:=1/2-\eps$ and
\begin{equation}\label{old-param}
p_0:={1\over10}-\eps,\qquad p_1=p_2:={3\over4}-\eps,
\end{equation}
where $\eps>0$ is a small bias parameter (less than 1/10, of course).  In game $A$, the player tosses a $p$-coin (i.e., $p$ is the probability of heads).  In game $B$, if the player's current capital is congruent to $j$ (mod 3), he tosses a $p_j$-coin.  (Assume initial capital 0 for simplicity.) In both games, the player wins one unit with heads and loses one unit with tails.

It can be shown that games $A$ and $B$ are both losing games (asymptotically), regardless of $\eps$, whereas the random mixture $(1/2)(A+B)$ (toss a fair coin to determine which game to play) is a winning game for $\eps$ sufficiently small.  Furthermore, certain nonrandom patterns, including $AAB$, $ABB$, and $AABB$ but excluding $AB$, are winning as well, again for $\eps$ sufficiently small.  These are the original examples of \textit{Parrondo's paradox}.

It has been suggested that game $A$ acts as ``noise'' to break up the losing cycles of game $B$ played alone (Kay and Johnson 2003).  Toral (2002) proposed a stochastic model in which a different type of noise appears to have a similar effect.  The model assumes an ensemble of $N\ge2$ players and replaces the noise effect of Parrondo's game $A$ by a redistribution of capital among the players.  A player $i$ is selected at random to play. With probability 1/2 he can either play Parrondo's game $B$ or game $A'$ consisting in that player giving away one unit of his capital to a randomly selected (without replacement) player $j$.  Notice that this new game $A'$ is fair since it does not modify the total amount of capital, it simply redistributes it randomly among the players.

Toral showed by computer simulation that the Parrondo effect is present in his games.  Our aim here is to prove this, establishing a strong law of large numbers and a central limit theorem for the sequence of profits of the ensemble of $N$ players.  For this we apply results of Ethier and Lee (2009), but the application is not straightforward.  For example, the formulas for the mean and variance parameters in the central limit theorem depend on the unique stationary distribution of the underlying Markov chain as well as on its fundamental matrix, both of which are too complicated to derive explicitly except for small $N$.  Nevertheless, we can evaluate the mean and variance parameters for all $N$.

We generalize (\ref{old-param}) to the parameterization of Ethier and Lee (2009):
\begin{equation}\label{param}
p_0:={\rho^2\over1+\rho^2}-\eps,\qquad p_1=p_2:={1\over1+\rho}-\eps,
\end{equation}
where $\rho>0$ (eq.\ (\ref{old-param}) is the special case $\rho=1/3$).  The bias parameter is not important, so we take $\eps=0$ in most of what follows, which makes game $B$ fair (asymptotically).

Let us summarize our results.  Just as it is conventional in the literature to denote the nonrandom pattern $(A')^r B^s$ by $[r,s]$, we will introduce the (slightly redundant) notation $(\gamma,1-\gamma)$ for the random mixture $\gamma A'+(1-\gamma)B$.  We establish a strong law of large numbers (SLLN) and a central limit theorem (CLT) for the sequence of profits of the ensemble of $N$ players in both settings (random mixture and nonrandom pattern).  We provide a formula for the random-mixture mean $\mu_{(\gamma,1-\gamma)}^{(N)}$, which does not depend on $N$, as a function of $\gamma\in(0,1)$ and $\rho>0$.  The nonrandom-pattern mean $\mu_{[r,s]}^{(N)}$ does depend on $N$ and is rather more complicated; we provide a formula, as a function of $N\ge2$ and $\rho>0$, only for small $r,s\ge1$ but we determine its sign for all $r,s\ge1$, $N\ge2$, and $\rho>0$, thereby establishing necessary and sufficient conditions for the Parrondo effect to be present.  Finally we show that the random-mixture case and the nonrandom-pattern case are connected by the unexpected relationship
\begin{equation}\label{relation}
\mu_{(r/(r+s),s/(r+s))}^{(N)}=\lim_{M\to\infty}\mu_{[r,s]}^{(M)},\qquad r,s\ge1,\;N\ge2,\;\rho>0,
\end{equation}
and a simple formula for this common value is provided.  To put this in perspective, the corresponding identity for one-player Parrondo games fails in all but one case ($r=2$, $s=1$).

The variance parameter is considerably more complicated, so we assume that $\rho=1/3$ (i.e., (\ref{old-param}) holds with $\eps=0$) and $\gamma=1/2$, obtaining a formula for $(\sigma_{(1/2,1/2)}^{(N)})^2$ as a function of $N\ge2$.  We do the same for $(\sigma_{[r,s]}^{(N)})^2$ for $\rho=1/3$ and small $r,s\ge1$.  It turns out that the analogue of (\ref{relation}) fails for the variances.  However, a different notion of variance, the expected sample variance of the individual players' capitals, which was considered by Toral (2002), does apparently satisfy a relationship nearly analogous to (\ref{relation}).  We can confirm this only in special cases, so it remains a conjecture.

Toral (2002) also studied a model in which the capital-dependent game is replaced by a history-dependent game.  It seems likely that most of the results of this paper can be extended to that setting, with the probable exception of Theorem \ref{positivity} below.  Finally, Toral considered a model with redistribution of wealth from richer to poorer neighbors.  That model is considerably more difficult to analyze than either of the others, and no rigorous results for it are known.

\section{Mean profit for random mixtures}\label{model}

There are two natural ways to define the model.  The simplest is to describe the state of the system by an $N$-dimensional vector $\bm x=(x_1,x_2,\ldots,x_N)$ in which $x_i$ denotes the capital (mod 3) of player $i$.  An alternative approach (adopted by Ethier 2007), which makes the state space smaller but the one-step transition probabilities more complicated, is to describe the state of the system, when it is in state $\bm x$ according to the previous description, by $(n_0,n_1,n_2)$, where $n_0$ (resp., $n_1$, $n_2$) is the number of 0s (resp., 1s, 2s) among $x_1,x_2,\ldots,x_N$.  Using the first approach, the state space is
$$
\Sigma_N:=\{\bm x=(x_1,x_2,\ldots,x_N): x_i\in\{0,1,2\}{\rm\ for\ }i=1,\ldots,N\}=\{0,1,2\}^N,
$$
while using the second approach, the state space is
$$
\bar\Sigma_N:=\{(n_0,n_1,n_2)\in{\bf Z}_+^3: n_0+n_1+n_2=N\}.
$$
We note that $|\Sigma_N|=3^N$ and $|\bar\Sigma_N|={N+2\choose2}$.

The one-step transition probabilities using the first approach depend on three probabilities $p_0,p_1,p_2$.  If only game $B$ is played, then the one-step transition probabilities have the simple form
$$
\bm P_B^{(N)}(\bm x,\bm y):=\begin{cases}N^{-1}p_{x_i}&\text{if $y_i=x_i+1$ (mod 3) and $y_j=x_j$ for all $j\ne i$}\\
N^{-1}q_{x_i}&\text{if $y_i=x_i-1$ (mod 3) and $y_j=x_j$ for all $j\ne i$}\end{cases}
$$
for $i=1,2,\ldots,N$, where $q_x:=1-p_x$ for $x=0,1,2$, and $\bm P_B^{(N)}(\bm x,\bm y)=0$ otherwise.  We adopt the parameterization (\ref{param}) with $\eps=0$.

If only game $A'$ is played, then the one-step transition matrix is symmetric and of the form
$$
\bm P_{A'}^{(N)}(\bm x,\bm y):=[N(N-1)]^{-1}
$$
if, for some $i,j\in\{1,2,\ldots,N\}$ with $i\ne j$, we have $y_i=x_i-1$ (mod 3), $y_j=x_j+1$ (mod 3), and $y_k=x_k$ for all $k\ne i,j$.  Finally, if the two games are mixed, that is, game $A'$ is played with probability $\gamma\in(0,1)$ and game $B$ is played with probability $1-\gamma$, then our one-step transition matrix has the form $\bm P_{(\gamma,1-\gamma)}^{(N)}:=\gamma\bm P_{A'}^{(N)}+(1-\gamma)\bm P_B^{(N)}$.

The one-step transition probabilities using the second approach also depend on the three probabilities $p_0,p_1,p_2$ and are best summarized in the form of a table.  See Table \ref{Table1}, which is  essentially from Ethier (2007).

\begin{table}[ht]
\caption{\label{Table1}One-step transitions using the second approach, for both game $A'$ and game $B$.  From state $(n_0,n_1,n_2)$, a transition is made to state $(n_0',n_1',n_2')$.\medskip}
\tabcolsep=.16cm
\begin{center}
\begin{tabular}{ccccc}
\hline
\noalign{\smallskip}
&&& type of &\\
$(n_0',n_1',n_2')$  & type of & game   & winner  & probability \\
      & player & played & / result &             \\
\noalign{\smallskip}
\hline
\noalign{\smallskip}
$(n_0-2,n_1+1,n_2+1)$ & 0 & $A'$ & 0 & $[N(N-1)]^{-1}n_0(n_0-1)$ \\
$(n_0-1,n_1-1,n_2+2)$ & 0 & $A'$ & 1 & $[N(N-1)]^{-1}n_0 n_1$ \\
$(n_0,n_1,n_2)$ & 0 & $A'$ & 2 & $[N(N-1)]^{-1}n_0 n_2$ \\
$(n_0,n_1,n_2)$ & 1 & $A'$ & 0 & $[N(N-1)]^{-1}n_1 n_0$ \\
$(n_0+1,n_1-2,n_2+1)$ & 1 & $A'$ & 1 & $[N(N-1)]^{-1}n_1(n_1-1)$ \\
$(n_0+2,n_1-1,n_2-1)$ & 1 & $A'$ & 2 & $[N(N-1)]^{-1}n_1 n_2$ \\
$(n_0-1,n_1+2,n_2-1)$ & 2 & $A'$ & 0 & $[N(N-1)]^{-1}n_2 n_0$ \\
$(n_0,n_1,n_2)$ & 2 & $A'$ & 1 & $[N(N-1)]^{-1}n_2 n_1$ \\
$(n_0+1,n_1+1,n_2-2)$ & 2 & $A'$ & 2 & $[N(N-1)]^{-1}n_2(n_2-1)$ \\
\noalign{\smallskip}
\hline
\noalign{\smallskip}
$(n_0-1,n_1+1,n_2)$ & 0 & $B$ & win & $N^{-1}n_0 p_0$ \\
$(n_0-1,n_1,n_2+1)$ & 0 & $B$ & lose & $N^{-1}n_0 q_0$ \\
$(n_0,n_1-1,n_2+1)$ & 1 & $B$ & win & $N^{-1}n_1 p_1$ \\
$(n_0+1,n_1-1,n_2)$ & 1 & $B$ & lose & $N^{-1}n_1 q_1$ \\
$(n_0+1,n_1,n_2-1)$ & 2 & $B$ & win & $N^{-1}n_2 p_2$ \\
$(n_0,n_1+1,n_2-1)$ & 2 & $B$ & lose & $N^{-1}n_2 q_2$ \\
\noalign{\smallskip}
\hline
\end{tabular}
\end{center}
\end{table}

That the two approaches to the model are equivalent, at least in the stationary setting, is a consequence of the following simple lemma, which is easily seen to be applicable to $\bm P_B^{(N)}$ and $\bm P_{(\gamma,1-\gamma)}^{(N)}$.

We first need some notation.  Given a finite set $E$ and an integer $N\ge2$, put $E^N:=E\times\cdots\times E$.  Given a permutation $\sigma$ of $\{1,2,\ldots,N\}$ and $\bm x=(x_1,\ldots,x_N)\in E^N$, write $\bm x_\sigma:=(x_{\sigma(1)},\ldots,x_{\sigma(N)})$.

\begin{lemma}
Let $E$ be a finite set, fix $N\ge2$, let $\bm P$ be the one-step transition matrix for an irreducible Markov chain in the product space $E^N$, and let $\bm\pi$ be its unique stationary distribution.  If, for every permutation $\sigma$ of $\{1,2,\ldots,N\}$,
$$
\bm P(\bm x_\sigma,\bm y_\sigma)=\bm P(\bm x,\bm y)
$$
for all $\bm x,\bm y\in E^N$, then $\bm\pi$ is exchangeable, that is, for every permutation $\sigma$ of $\{1,2,\ldots,N\}$, we have $\bm\pi(\bm x_\sigma)=\bm\pi(\bm x)$ for all $\bm x\in E^N$.
\end{lemma}

\begin{proof}  Given a permutation $\sigma$ of $\{1,2,\ldots,N\}$, define the distribution ${\bm\pi}_\sigma$ on $E^N$ by ${\bm\pi}_\sigma(\bm x):=\bm\pi(\bm x_\sigma)$.  Then
\begin{eqnarray*}
\bm\pi_\sigma(\bm y)
=\sum_{\bm x\in E^N}\bm\pi(\bm x)\bm P(\bm x,\bm y_\sigma)=\sum_{\bm x\in E^N}\bm\pi(\bm x_\sigma)\bm P(\bm x_\sigma,\bm y_\sigma)
=\sum_{\bm x\in E^N}\bm\pi_\sigma(\bm x)\bm P(\bm x,\bm y)
\end{eqnarray*}
for all $\bm y\in E^N$, hence by the uniqueness of stationary distributions, $\bm\pi_\sigma=\bm\pi$.
\end{proof}

We would like to apply results of Ethier and Lee (2009) to game $B$ and to the mixed game.  (They do not apply to game $A'$ because the one-step transition matrix $\bm P_{A'}^{(N)}$ is not irreducible, but the behavior of the system is clear in this case.)  We restate those results here for convenience.

Consider an irreducible aperiodic Markov chain $\{X_n\}_{n\ge0}$ with finite state space $\Sigma$.  It evolves according to the one-step transition matrix ${\bm P}=(P_{ij})_{i,j\in\Sigma}$.  Let us denote its unique stationary distribution by ${\bm \pi}=(\pi_i)_{i\in \Sigma}$.  Let $w:\Sigma\times\Sigma\mapsto {\bf R}$ be an arbitrary function, which we write as a matrix ${\bm W}=(w(i,j))_{i,j\in\Sigma}$ and refer to as the \textit{payoff matrix}. Finally, define the sequences $\{\xi_n\}_{n\ge1}$ and $\{S_n\}_{n\ge1}$ by
\begin{equation}\label{xi_n}
\xi_n:=w(X_{n-1},X_n),\qquad n\ge1,
\end{equation}
and
\begin{equation}\label{S_n}
S_n:=\xi_1+\cdots+\xi_n,\qquad n\ge1.
\end{equation}
Let ${\bm \Pi}$ denote the square matrix each of whose rows is ${\bm \pi}$, and let ${\bm Z}:=({\bm I}-({\bm P}-{\bm \Pi}))^{-1}$ denote the \textit{fundamental matrix}.  Denote by $\dot{\bm P}$ (resp., $\ddot{\bm P}$) the Hadamard (entrywise) product $\bm P\circ\bm W$ (resp., $\bm P\circ\bm W\circ\bm W$), and let $\bm 1:=(1,1,\ldots,1)^\T$.  Then define
\begin{equation}\label{mu,sigma2}
\mu:=\bm\pi\dot{\bm P}\bm 1\quad{\rm and}\quad\sigma^2:=\bm\pi\ddot{\bm P}\bm 1
-(\bm\pi\dot{\bm P}\bm 1)^2+2\bm\pi\dot{\bm P}(\bm Z-\bm\Pi)\dot{\bm P}\bm 1.
\end{equation}

\begin{theorem}[Ethier and Lee 2009]\label{EL}
Under the above assumptions, and with the distribution of $X_0$ arbitrary, $\lim_{n\to\infty}n^{-1}\E[S_n]=\mu$,
$$
{S_n\over n}\to \mu\;\;{\rm a.s.},
$$
$\lim_{n\to\infty}n^{-1}\Var(S_n)=\sigma^2$, and, if $\sigma^2>0$,
$$
{S_n-n\mu\over\sqrt{n\sigma^2}}\to_d N(0,1).
$$
If $\mu=0$ and $\sigma^2>0$, then $-\infty=\liminf_{n\to\infty}S_n<\limsup_{n\to\infty}S_n=\infty$ \emph{a.s.}
\end{theorem}

We apply this result first with $\Sigma:=\Sigma_N$ and $\bm P:=\bm P_B^{(N)}$, which is clearly irreducible and aperiodic.  We claim that the stationary distribution $\bm\pi_B^{(N)}$ is the $N$-fold product measure $\bm\pi\times\bm\pi\times\cdots\times\bm\pi$, where $\bm\pi=(\pi_0,\pi_1,\pi_2)$ denotes the stationary distribution of the three-state chain in $\Sigma_1$ with one-step transition matrix
$$
\bm P_B^{(1)}=\left(\begin{array}{ccc}
0&p_0&q_0\\
q_1&0&p_1\\
p_2&q_2&0
\end{array}\right).
$$
Indeed,
\begin{eqnarray*}
&&\sum_{\bm x}\pi_{x_1}\cdots\pi_{x_N}\bm P_B^{(N)}(\bm x,\bm y)\\
&&\quad{}=\sum_{i=1}^N\pi_{y_1}\cdots\pi_{y_{i-1}}\pi_{y_{i+1}}\cdots\pi_{y_N}\\
\noalign{\vglue-3mm}
&&\qquad\qquad\qquad\qquad\qquad{}\cdot\sum_{x_i:x_i\ne y_i}\pi_{x_i}\bm P_B^{(N)}((y_1,\ldots,y_{i-1},x_i,y_{i+1},\ldots,y_N),\bm y)\\
&&\quad{}=N^{-1}\sum_{i=1}^N\pi_{y_1}\cdots\pi_{y_N}\\
&&\quad{}=\pi_{y_1}\cdots\pi_{y_N},
\end{eqnarray*}
where the first equality holds because state $\bm y$ can be reached in one step only from states $\bm x$ that differ from $\bm y$ at exactly one coordinate.
Alternatively, we could take $\Sigma:=\bar\Sigma_N$ and $\bm P:=\bar{\bm P}_B^{(N)}$ from Table \ref{Table1}.  In this case the unique stationary distribution is multinomial$(N,\bm\pi)$.

Next, let us determine the value of $\mu$ in the theorem.  We have
\begin{eqnarray*}
\mu_B^{(N)}&=&\bm\pi_B^{(N)}\dot{\bm P}_B^{(N)}\bm 1=\sum_{\bm x}\pi_{x_1}\cdots\pi_{x_N}\sum_{i=1}^N N^{-1}(p_{x_i}-q_{x_i})\\
&=&N^{-1}\sum_{(n_0,n_1,n_2)}{N\choose n_0,n_1,n_2}\pi_0^{n_0}\pi_1^{n_1}\pi_2^{n_2}[n_0(p_0-q_0)+n_1(p_1-q_1)\\
\noalign{\vglue-4mm}
&&\qquad\qquad\qquad\qquad\qquad\qquad\qquad\qquad\qquad\qquad\qquad\qquad{}+n_2(p_2-q_2)]\\
&=&\pi_0(p_0-q_0)+\pi_1(p_1-q_1)+\pi_2(p_2-q_2)=\mu_B^{(1)}=0
\end{eqnarray*}
because the parameterization (\ref{param}) with $\eps=0$ was chosen to ensure the last equality.

Now we apply the theorem with $\Sigma:=\Sigma_N$ and $\bm P:=\bm P_{(\gamma,1-\gamma)}^{(N)}=\gamma\bm P_{A'}^{(N)}+(1-\gamma)\bm P_B^{(N)}$, where $0<\gamma<1$, which is also irreducible and aperiodic (because $\bm P_B^{(N)}$ is).  Here the unique stationary distribution $\bm\pi_{(\gamma,1-\gamma)}^{(N)}$ is complicated.  For example, in the simplest case, $\gamma=1/2$ and $N=2$,
\begin{eqnarray*}
\bm\pi_{(1/2,1/2)}^{(2)}(0,0)&=&(1 + \rho^2) (31 + 47 \rho + 60 \rho^2 + 47 \rho^3 + 31 \rho^4)/d,\\
\bm\pi_{(1/2,1/2)}^{(2)}(0,1)&=&\bm\pi_{(1/2,1/2)}^{(2)}(1,0)=2(1 + \rho) (1 + \rho^2) (11 + 15 \rho + 9 \rho^2 + 19 \rho^3)/d,\\
\bm\pi_{(1/2,1/2)}^{(2)}(0,2)&=&\bm\pi_{(1/2,1/2)}^{(2)}(2,0)=2(1 + \rho) (1 + \rho^2) (19 + 9 \rho + 15 \rho^2 + 11 \rho^3)/d,\\
\bm\pi_{(1/2,1/2)}^{(2)}(1,1)&=&(1 + \rho) (19 + 21 \rho + 48 \rho^2 + 59 \rho^3 + 27 \rho^4 + 42 \rho^5)/d,\\
\bm\pi_{(1/2,1/2)}^{(2)}(1,2)&=&\bm\pi_{(1/2,1/2)}^{(2)}(2,1)=6 (1 + \rho)^2 (1 + \rho^2) (4 + \rho + 4 \rho^2)/d,\\
\bm\pi_{(1/2,1/2)}^{(2)}(2,2)&=&(1 + \rho) (42 + 27 \rho + 59 \rho^2 + 48 \rho^3 + 21 \rho^4 + 19 \rho^5)/d,
\end{eqnarray*}
where $d:=2 (13 - 2 \rho + 13 \rho^2) (10 + 20 \rho + 21 \rho^2 + 20 \rho^3 + 10 \rho^4)$.  In particular, each entry of $\bm\pi_{(1/2,1/2)}^{(2)}$ is the ratio of two degree-6 polynomials in $\rho$.  In another simple case, $\gamma=1/2$ and $N=3$, each entry of $\bm\pi_{(1/2,1/2)}^{(3)}$ is the ratio of two degree-14 polynomials in $\rho$.  Fortunately, explicit formulas such as these are unnecessary to evaluate $\mu_{(\gamma,1-\gamma)}^{(N)}$.

Let $\bar{\bm\pi}_{(\gamma,1-\gamma)}^{(N)}$ denote the corresponding stationary distribution on $\bar\Sigma_N$.  Then the mean profit per turn to the ensemble of players is
\begin{eqnarray}\label{mu_C-prelim}
\mu_{(\gamma,1-\gamma)}^{(N)}&=&\bm\pi_{(\gamma,1-\gamma)}^{(N)}\dot{\bm P}_{(\gamma,1-\gamma)}^{(N)}\bm 1\nonumber\\
&=&(1-\gamma)\sum_{\bm x}\bm\pi_{(\gamma,1-\gamma)}^{(N)}(x_1,\ldots,x_N)\sum_{i=1}^N N^{-1}(p_{x_i}-q_{x_i})\\
&=&N^{-1}(1-\gamma)\sum_{(n_0,n_1,n_2)}\bar{\bm\pi}_{(\gamma,1-\gamma)}^{(N)}(n_0,n_1,n_2) [n_0(p_0-q_0)+n_1(p_1-q_1)\nonumber\\
\noalign{\vglue-4mm}
&&\qquad\qquad\qquad\qquad\qquad\qquad\qquad\qquad\qquad\qquad\qquad\;\;{}+n_2(p_2-q_2)]\nonumber\\
&=&N^{-1}(1-\gamma)\{\E_{\bar{\bm\pi}_{(\gamma,1-\gamma)}^{(N)}}[n_0](p_0-q_0)+\E_{\bar{\bm\pi}_{(\gamma,1-\gamma)}^{(N)}}[n_1](p_1-q_1)\nonumber\\
&&\qquad\qquad\qquad\qquad\qquad\qquad\qquad\qquad\qquad{}+\E_{\bar{\bm\pi}_{(\gamma,1-\gamma)}^{(N)}}[n_2](p_2-q_2)\}.\nonumber
\end{eqnarray}
Now by Table \ref{Table1}, we can compute
\begin{eqnarray*}
\E[n_0'-n_0]&=&\gamma{-2n_0(n_0-1)-n_0n_1+n_1(n_1-1)+2n_1n_2-n_2n_0+n_2(n_2-1)\over N(N-1)}\\
&&\quad{}+(1-\gamma){-n_0p_0-n_0q_0+n_1q_1+n_2p_2\over N}\\
&=&{\gamma(N-3n_0)+(1-\gamma)[n_0(-1)+n_1q_1+n_2p_2]\over N}.
\end{eqnarray*}
Similarly,
\begin{eqnarray*}
\E[n_1'-n_1]&=&{\gamma(N-3n_1)+(1-\gamma)[n_0p_0+n_1(-1)+n_2q_2]\over N},\\
\E[n_2'-n_2]&=&{\gamma(N-3n_2)+(1-\gamma)[n_0q_0+n_1p_1+n_2(-1)]\over N}.
\end{eqnarray*}
In each of these equations, we have used $n_0+n_1+n_2=N$ to simplify, with the result that all the quadratic terms cancel and the right sides are linear in $(n_0,n_1,n_2)$, at least if we replace the $N$ in the numerators by $n_0+n_1+n_2$.

Next we take expectations with respect to $\bar{\bm\pi}_{(\gamma,1-\gamma)}^{(N)}$ to obtain
\begin{eqnarray*}
(0,0,0)&=&(\E_{\bar{\bm\pi}_{(\gamma,1-\gamma)}^{(N)}}[n_0],\E_{\bar{\bm\pi}_{(\gamma,1-\gamma)}^{(N)}}[n_1],\E_{\bar{\bm\pi}_{(\gamma,1-\gamma)}^{(N)}}[n_2])\left[\gamma\left(\begin{array}{rrr}
-2&1&1\\
1&-2&1\\
1&1&-2
\end{array}\right)\right.\\
&&\qquad\qquad\qquad\qquad\qquad\qquad\qquad\qquad\quad{}\left.{}+(1-\gamma)\left(\begin{array}{rrr}
-1&p_0&q_0\\
q_1&-1&p_1\\
p_2&q_2&-1
\end{array}\right)\right],
\end{eqnarray*}
which with $\E_{\bar{\bm\pi}_{(\gamma,1-\gamma)}^{(N)}}[n_0]+\E_{\bar{\bm\pi}_{(\gamma,1-\gamma)}^{(N)}}[n_1]+\E_{\bar{\bm\pi}_{(\gamma,1-\gamma)}^{(N)}}[n_2]=N$ uniquely determines the vector $(\E_{\bar{\bm\pi}_{(\gamma,1-\gamma)}^{(N)}}[n_0], \E_{\bar{\bm\pi}_{(\gamma,1-\gamma)}^{(N)}}[n_1],\E_{\bar{\bm\pi}_{(\gamma,1-\gamma)}^{(N)}}[n_2])$ because the matrix within brackets is an irreducible infinitesimal matrix.  Substituting into (\ref{mu_C-prelim}) and using our parametrization (\ref{param}) with $\eps=0$, we obtain
\begin{equation}\label{mu_C}
\mu_{(\gamma,1-\gamma)}^{(N)}={3 \gamma(1-\gamma) (1 - \rho)^3 (1 + \rho)\over
2 (1 + \rho + \rho^2)^2 + \gamma (5 + 10 \rho + 6 \rho^2 + 10 \rho^3 + 5 \rho^4) + 2 \gamma^2 (1 + \rho + \rho^2)^2},
\end{equation}
which does not depend on $N$ and is positive if $0<\rho<1$, zero if $\rho=1$, and negative if $\rho>1$, indicating that the Parrondo effect is present, regardless of $\gamma\in(0,1)$, if $\rho\ne1$.  (In the case $\rho>1$, the effect is sometimes referred to as a \textit{reverse} Parrondo effect.  We will not make this distinction.)  Temporarily denoting $\mu_{(\gamma,1-\gamma)}^{(N)}$ by $\mu_{(\gamma,1-\gamma)}^{(N)}(\rho)$ to emphasize its dependence on $\rho$, we note that
$$
\mu_{(\gamma,1-\gamma)}^{(N)}(1/\rho)=-\mu_{(\gamma,1-\gamma)}^{(N)}(\rho),
$$
a fact that can also be proved probabilistically (Ethier and Lee 2009).

When $\gamma=1/2$, this reduces to
$$
\mu_{(1/2,1/2)}^{(N)}={3 (1 - \rho)^3 (1 + \rho)\over 2(10 + 20 \rho + 21 \rho^2 + 20 \rho^3 + 10 \rho^4)}.
$$
As we will see in Section \ref{coincidence}, this formula appears elsewhere in the literature of Parrondo's paradox.

\section{An alternative approach}\label{alternative}

The method used in Section \ref{model} to find $\mu_{(\gamma,1-\gamma)}^{(N)}$ does not extend to finding the variance $(\sigma_{(\gamma,1-\gamma)}^{(N)})^2$.  However, a method that does extend is based on the observation that the components of the $N$-dimensional Markov chain controlling the mixed game are themselves Markovian.

For example, when game $B$ is played, the Markov chain for player $i$ (one of the $N$ players) has one-step transition matrix
\begin{equation}\label{pB1N}
\bm P_B^{(1,N)}:=N^{-1}[\bm P_B^{(1)}+(N-1)\bm I_3].
\end{equation}
On the other hand, the redistribution game $A'$ affects player $i$ only if $i$ is chosen as the donor or as the beneficiary (probability $(N-1)/[N(N-1)]=1/N$ for each).  This leads to
\begin{equation}\label{pA1N}
\bm P_{A'}^{(1,N)}:=N^{-1}[2\bm P_A^{(1)}+(N-2)\bm I_3],
\end{equation}
where $\bm P_A^{(1)}$ denotes the one-step transition matrix for the original one-player Parrondo game $A$ (not $A'$).  In both displayed matrices, the superscript $(1,N)$ is intended to indicate that the underlying Markov chain controls one of the $N$ players.

From these one-step transition matrices we calculate
$$
\dot{\bm P}_B^{(1,N)}:=N^{-1}\dot{\bm P}_B^{(1)},
\qquad
\dot{\bm P}_{A'}^{(1,N)}:=2N^{-1}\dot{\bm P}_A^{(1)},
$$
and
$$
\ddot{\bm P}_B^{(1,N)}:=N^{-1}\ddot{\bm P}_B^{(1)},
\qquad
\ddot{\bm P}_{A'}^{(1,N)}:=2N^{-1}\ddot{\bm P}_A^{(1)}.
$$
With
\begin{eqnarray*}
\bm P&:=&\gamma\bm P_{A'}^{(1,N)}+(1-\gamma)\bm P_B^{(1,N)},\\
\dot{\bm P}&:=&\gamma\dot{\bm P}_{A'}^{(1,N)}+(1-\gamma)\dot{\bm P}_B^{(1,N)},\\
\ddot{\bm P}&:=&\gamma\ddot{\bm P}_{A'}^{(1,N)}+(1-\gamma)\ddot{\bm P}_B^{(1,N)},
\end{eqnarray*}
and with $\bm\pi$, $\bm\Pi$, and $\bm Z$ chosen accordingly and $\bm 1:=(1,1,1)^\T$, we have
$$
\mu_{(\gamma,1-\gamma)}^{(1,N)}=\bm\pi\dot{\bm P}\bm1,\qquad (\sigma_{(\gamma,1-\gamma)}^{(1,N)})^2=\bm\pi\ddot{\bm P}\bm1-(\bm\pi\dot{\bm P}\bm1)^2+2\bm\pi\dot{\bm P}(\bm Z-\bm\Pi)\dot{\bm P}\bm1.
$$
The mean is readily evaluated to give
\begin{eqnarray}\label{mu_C-2}
\mu_{(\gamma,1-\gamma)}^{(N)}&=&N\mu_{(\gamma,1-\gamma)}^{(1,N)}\\
&=&{3 \gamma(1 - \gamma)  (1 - \rho)^3 (1 + \rho)\over 2 (1 + \rho + \rho^2)^2 +
   \gamma (5 + 10 \rho + 6 \rho^2 + 10 \rho^3 + 5 \rho^4)+ 2 \gamma^2 (1 + \rho + \rho^2)^2},\nonumber
\end{eqnarray}
which is consistent with (\ref{mu_C}) and does not depend on $N$.  The variance $(\sigma_{(\gamma,1-\gamma)}^{(1,N)})^2$ is also easily evaluated but is  complicated; instead we provide its asymptotic value as $N\to\infty$ ($a_N\sim b_N$ if $\lim_{N\to\infty}a_N/b_N=1$):
\begin{eqnarray}\label{sigma_C-asymp}
&&\!\!\!\!\!(\sigma_{(\gamma,1-\gamma)}^{(1,N)})^2\nonumber\\
&&\!\!\!\!\!{}\sim 9 [8 (1 + \gamma^7) \rho^2 (1 + \rho + \rho^2)^4 \nonumber\\
&&\;\;{} + 4 (\gamma + \gamma^6) (1 + \rho + \rho^2)^2 (1 + 2 \rho + \rho^2 + 2 \rho^3 + \rho^4) (1 + 2 \rho + 12 \rho^2 + 2 \rho^3 + \rho^4)\nonumber\\
&&\;\;{} + 6 (\gamma^2 + \gamma^5) (1 + \rho + \rho^2)^2 (3 + 20 \rho + 30 \rho^2 + 40 \rho^3 + 66 \rho^4 + 40 \rho^5 + 30 \rho^6 \nonumber\\
&&\qquad\qquad\qquad\qquad\qquad\qquad\;\;{}+ 20 \rho^7 + 3 \rho^8) \nonumber\\
&&\;\;{} + (\gamma^3 + \gamma^4) (59 + 306 \rho + 864 \rho^2 + 1738 \rho^3 + 2781 \rho^4 + 3636 \rho^5 + 3912 \rho^6 \nonumber\\
&&\qquad\qquad\qquad\;{} + 3636 \rho^7 + 2781 \rho^8 + 1738 \rho^9 + 864 \rho^{10} + 306 \rho^{11} + 59 \rho^{12})]\nonumber\\
&&\;{} /\{N [2 (1 + \gamma^2 ) (1 + \rho + \rho^2)^2 + \gamma (5 + 10 \rho + 6 \rho^2 + 10 \rho^3 + 5 \rho^4)]^3\}.
\end{eqnarray}

\section{Variance parameter for game $B$}\label{varB}

Let $\bm P$ be the one-step transition matrix for an irreducible aperiodic Markov chain, let $\bm\pi$ be its unique stationary distribution, and let $\bm\Pi$ be the square matrix each of whose rows is $\bm\pi$.  Denote by $\bm Z_{\bm P}:=(\bm I-(\bm P-\bm\Pi))^{-1}$ the fundamental matrix of $\bm P$.

\begin{lemma}\label{fundamental}For each positive integer $N$,
$\bm Z_{(1/N)\bm P+(1-1/N)\bm I}-\bm\Pi=N(\bm Z_{\bm P}-\bm\Pi)$.
\end{lemma}

\begin{proof}
The one-step transition matrix $(1/N)\bm P+(1-1/N)\bm I$ has the same stationary distribution $\bm\pi$, hence the same $\bm\Pi$, so
$$
\bm Z_{(1/N)\bm P+(1-1/N)\bm I}=(\bm I-[(1/N)\bm P+(1-1/N)\bm I-\bm\Pi])^{-1}=N(\bm I-(\bm P-N\bm\Pi))^{-1},
$$
hence it suffices to prove that
$$
(\bm I-(\bm P-N\bm\Pi))^{-1}-(1/N)\bm\Pi=(\bm I-(\bm P-\bm\Pi))^{-1}-\bm\Pi.
$$
For this it is enough that
\begin{eqnarray*}
&&(\bm I-(\bm P-N\bm\Pi))[(\bm I-(\bm P-N\bm\Pi))^{-1}-(1/N)\bm\Pi]\\
&&\qquad{}=(\bm I-(\bm P-\bm\Pi)+(N-1)\bm\Pi)[(\bm I-(\bm P-\bm\Pi))^{-1}-\bm\Pi]
\end{eqnarray*}
or
\begin{eqnarray}\label{eq}
&&\bm I-(1/N)(\bm I-(\bm P-N\bm\Pi))\bm\Pi\nonumber\\
&&\qquad{}=\bm I-(\bm I-(\bm P-\bm\Pi))\bm\Pi+(N-1)\bm\Pi[(\bm I-(\bm P-\bm\Pi))^{-1}-\bm\Pi].
\end{eqnarray}
Now $\bm\Pi\bm P=\bm P\bm\Pi=\bm\Pi$, $\bm\Pi^2=\bm\Pi$, and so $\bm\Pi=\bm\Pi(\bm I-(\bm P-\bm\Pi))$ and $\bm\Pi(\bm I-(\bm P-\bm\Pi))^{-1}=\bm\Pi$.  So (\ref{eq}) is equivalent to
$$
\bm I-(1/N)(\bm\Pi-(\bm\Pi-N\bm\Pi))=\bm I-(\bm\Pi-(\bm\Pi-\bm\Pi))+(N-1)(\bm\Pi-\bm\Pi)
$$
or $\bm I-\bm\Pi=\bm I-\bm\Pi$, hence (\ref{eq}), and therefore the lemma, is established.
\end{proof}

We want to use this to evaluate the variance parameter for Toral's $N$-player game $B$, in which there is no redistribution of wealth.  The state space is $\Sigma_N$ and the one-step transition probabilities are as previously specified.  We assume the parameterization (\ref{param}) with $\eps=0$.

We have seen that the stationary distribution $\bm\pi_B^{(N)}$ is the $N$-fold product measure $\bm\pi\times\bm\pi\times\cdots\times\bm\pi$, where $\bm\pi=(\pi_0,\pi_1,\pi_2)$ denotes the stationary distribution of the three-state chain with one-step transition matrix $\bm P_B^{(1)}$.
Specifically,
$$
\pi_0={1+\rho^2\over2(1+\rho+\rho^2)},\quad \pi_1={\rho(1+\rho)\over2(1+\rho+\rho^2)},\quad \pi_2={1+\rho\over2(1+\rho+\rho^2)}.
$$
In principle, we could use the formula
$\sigma^2:=\bm\pi\ddot{\bm P}\bm 1-(\bm\pi\dot{\bm P}\bm 1)^2+2\bm\pi\dot{\bm P}(\bm Z-\bm\Pi)\dot{\bm P}\bm 1$,
but the evaluation of the $3^N\times 3^N$ fundamental matrix $\bm Z$ is difficult, so we take a different approach.

The key observation is that each coordinate of the $N$-dimensional Markov chain is a one-dimensional Markov chain with one-step transition matrix (\ref{pB1N}) or
$$
\bm P_B^{(1,N)}:=(1/N)\bm P_B^{(1)}+(1-1/N)\bm I_3.
$$
Further, the coordinate processes are independent if their initial states are, and they are if the initial state of the $N$-dimensional process has the  stationary distribution $\bm\pi_B^{(N)}$ on $\Sigma_N$.

As already noted in Section \ref{alternative}, $\dot{\bm P}_B^{(1,N)}=(1/N)\dot{\bm P}_B^{(1)}$ and $\ddot{\bm P}_B^{(1,N)}=(1/N)\ddot{\bm P}_B^{(1)}$.  By Lemma \ref{fundamental}, $\bm Z_B^{(1,N)}-\bm\Pi=N(\bm Z_B^{(1)}-\bm\Pi)$,
so (since $\mu_B^{(1,N)}=N^{-1}\mu_B^{(1)}=0$)
\begin{eqnarray*}
(\sigma_B^{(1,N)})^2&:=&\bm\pi\ddot{\bm P}_B^{(1,N)}\bm 1+2\bm\pi\dot{\bm P}_B^{(1,N)}(\bm Z_B^{(1,N)}-\bm\Pi)\dot{\bm P}_B^{(1,N)}\bm 1\\
&\;=&N^{-1}[\bm\pi\ddot{\bm P}_B^{(1)}\bm 1+2\bm\pi\dot{\bm P}_B^{(1)}(\bm Z_B^{(1)}-\bm\Pi)\dot{\bm P}_B^{(1)}\bm 1]\\
&\;=&N^{-1}(\sigma_B^{(1)})^2=N^{-1}\bigg({3\rho\over1+\rho+\rho^2}\bigg)^2.
\end{eqnarray*}
Finally, let $S_n$ denote the profit to the ensemble of $N$ players after $n$ plays of game $B$, with $S_n^{[i]}$ denoting the profit to player $i$.  Then $S_n=S_n^{[1]}+\cdots+S_n^{[N]}$ and the summands are independent (assuming the stationary initial distribution mentioned above), hence
\begin{eqnarray}\label{sigmaB2}
(\sigma_B^{(N)})^2&=&\lim_{n\to\infty}n^{-1}\Var(S_n)=N\lim_{n\to\infty}n^{-1}\Var(S_n^{[1]})\nonumber\\
&=&N(\sigma_B^{(1,N)})^2=\bigg({3\rho\over1+\rho+\rho^2}\bigg)^2,
\end{eqnarray}
yielding a simple and explicit formula for $(\sigma_B^{(N)})^2$, which does not depend on $N$.

\section{Variance parameter for random mixtures}\label{variance_C}

With $S_n$ denoting the profit to the ensemble of $N$ players after $n$ plays of the mixed game, let $S_n^{[i]}$ denote the profit to player $i$ (one of the $N$ players) after $n$ plays of the mixed game.  Then
$$
S_n=\sum_{i=1}^N S_n^{[i]},
$$
so
\begin{eqnarray*}
\Var(S_n)&=&\sum_{i=1}^N\Var(S_n^{[i]})+2\sum_{1\le i<j\le N}\Cov(S_n^{[i]},S_n^{[j]})\\
&=&N\Var(S_n^{[1]})+N(N-1)\Cov(S_n^{[1]},S_n^{[2]}).
\end{eqnarray*}
Dividing by $n$ and letting $n\to\infty$, we find that
\begin{equation}\label{sigma_C^N}
(\sigma_{(\gamma,1-\gamma)}^{(N)})^2=N(\sigma_{(\gamma,1-\gamma)}^{(1,N)})^2+N(N-1)\sigma_{(\gamma,1-\gamma)}^{([1,2],N)},
\end{equation}
where the last superscript is intended to indicate that the underlying Markov chain controls players 1 and 2 of the $N$ players.
We know how to evaluate $(\sigma_{(\gamma,1-\gamma)}^{(1,N)})^2$, so it remains to find $\sigma_{(\gamma,1-\gamma)}^{([1,2],N)}$.

For this we will need an extension of (\ref{xi_n})--(\ref{mu,sigma2}).  With the same assumptions on $\{X_n\}_{n\ge0}$ (an irreducible, aperiodic, finite Markov chain in $\Sigma$ with one-step transition matrix $\bm P$ and unique stationary distribution $\bm\pi$), we let $w^{[1]},w^{[2]}:\Sigma\times\Sigma\mapsto{\bm R}$ be two functions with $\bm W^{[1]}$ and $\bm W^{[2]}$ denoting the corresponding matrices, and define
$$
\xi^{[1]}_n:=w^{[1]}(X_{n-1},X_n),\quad \xi^{[2]}_n:=w^{[2]}(X_{n-1},X_n),\qquad n\ge1,
$$
and
$$
S_n^{[1]}:=\xi^{[1]}_1+\cdots+\xi^{[1]}_n,\quad S_n^{[2]}:=\xi^{[2]}_1+\cdots+\xi^{[2]}_n,\qquad n\ge1.
$$
Let ${\bm\Pi}$ and ${\bm Z}$ be associated with $\bm P$ in the usual way.  Denote by $\bm P^{[1]}$, $\bm P^{[2]}$, and $\bm P^{[1,2]}$ the Hadamard products $\bm P\circ\bm W^{[1]}$, $\bm P\circ\bm W^{[2]}$, and $\bm P\circ\bm W^{[1]}\circ\bm W^{[2]}$, resp., and let $\bm 1:=(1,1,\ldots,1)^\T$.  Then define the covariance parameter
\begin{eqnarray*}
\sigma^{[1,2]}&:=&\bm\pi\bm P^{[1,2]}\bm 1-(\bm\pi\bm P^{[1]}\bm 1)(\bm\pi\bm P^{[2]}\bm 1)\nonumber\\
&&\quad{}+\bm\pi\bm P^{[1]}(\bm Z-\bm\Pi)\bm P^{[2]}\bm 1+\bm\pi\bm P^{[2]}(\bm Z-\bm\Pi)\bm P^{[1]}\bm 1.
\end{eqnarray*}
The interpretation of this parameter is as follows.

\begin{theorem}
Under the above assumptions, and with the distribution of $X_0$ arbitrary,
$$
\lim_{n\to\infty}n^{-1}\Cov(S^{[1]}_n,S^{[2]}_n)=\sigma^{[1,2]}.
$$
\end{theorem}

\begin{proof}
The proof is similar to the proof that $\lim_{n\to\infty}n^{-1}\Var(S_n)=\sigma^2$ in Theorem \ref{EL}, which is just the special case $w^{[1]}=w^{[2]}=w$.
\end{proof}

We now want to apply this to find $\sigma_{(\gamma,1-\gamma)}^{([1,2],N)}$.  This involves only players 1 and 2, for which we need only a (9-state) Markov chain in $\Sigma_2$.   The reduced model that does not distinguish between the players but only counts how many players of each type there are is insufficient.

Thinking of $(i,j)\in\Sigma_2$ as the base-3 representation of the integer $3i+j$, we order the elements of $\Sigma_2$ by their values (0--8).  The one-step transition matrix for the profit to players 1 and 2 when $N$ players are playing game $B$ is
$$
\bm P_B^{(2,N)}:=N^{-1}[2\bm P_B^{(2)}+(N-2)\bm I_9],
$$
where $\bm P_B^{(2)}$ is as in Section \ref{model} with $N=2$.  The superscript $(2,N)$ is intended to indicate that the underlying Markov chain controls two of the $N$ players.  The one-step transition matrix for the profit to players 1 and 2 when $N$ players are playing game $A'$ is
$$
\bm P_{A'}^{(2,N)}:=[N(N-1)]^{-1}[2\bm P_{A_0}+4(N-2)\bm P_{A_1}+(N-2)(N-3)\bm I_9],
$$
where $\bm P_{A_0}$ is a $9\times9$ matrix with two entries (each equal to 1/2) in each row, corresponding to one-unit transfers $1\to2$ and $2\to1$; similarly, $\bm P_{A_1}$ is a $9\times9$ matrix with four entries (each equal to 1/4) in each row, corresponding to one-unit transfers $1\to\cdot$, $\cdot\to1$, $2\to\cdot$, and $\cdot\to2$, where $\cdot$ represents the players other than 1 and 2.  The functions $w^{[1]}$ and $w^{[2]}$ can be specified as follows.  Corresponding to matrices $\bm P_B^{(2)}$ and $\bm P_{A_1}$, the function $w^{[1]}$ is 1 at (1 wins) and at $\cdot\to1$; it is $-1$ at (1 loses) and at $1\to\cdot$; and it is 0 at (2 wins) or (2 loses) and at $\cdot\to2$ and $2\to\cdot$.  Corresponding to matrix $\bm P_{A_0}$, the function $w^{[1]}$ is 1 at $2\to1$; it is $-1$ at $1\to2$.  The function $w^{[2]}$ is defined exactly in the same way but with the roles of 1 and 2 reversed.

From these one-step transition matrices we calculate
$$
({\bm P}_B^{(2,N)})^{[1]}:=2N^{-1}({\bm P}_B^{(2)})^{[1]},\qquad({\bm P}_B^{(2,N)})^{[2]}:=2N^{-1}({\bm P}_B^{(2)})^{[2]},
$$
\begin{eqnarray*}
({\bm P}_{A'}^{(2,N)})^{[1]}&:=&[N(N-1)]^{-1}[2({\bm P}_{A_0})^{[1]}+4(N-2)({\bm P}_{A_1})^{[1]}],\\
({\bm P}_{A'}^{(2,N)})^{[2]}&:=&[N(N-1)]^{-1}[2({\bm P}_{A_0})^{[2]}+4(N-2)({\bm P}_{A_1})^{[2]}],
\end{eqnarray*}
$({\bm P}_B^{(2,N)})^{[1,2]}:=\bm0$, and
$$
({\bm P}_{A'}^{(2,N)})^{[1,2]}:=2[N(N-1)]^{-1}({\bm P}_{A_0})^{[1,2]}.
$$
With
\begin{eqnarray*}
\bm P&:=&\gamma\bm P_{A'}^{(2,N)}+(1-\gamma)\bm P_B^{(2,N)},\\
{\bm P}^{[1]}&:=&\gamma({\bm P}_{A'}^{(2,N)})^{[1]}+(1-\gamma)({\bm P}_B^{(2,N)})^{[1]},\\
{\bm P}^{[2]}&:=&\gamma({\bm P}_{A'}^{(2,N)})^{[2]}+(1-\gamma)({\bm P}_B^{(2,N)})^{[2]},\\
{\bm P}^{[1,2]}&:=&\gamma({\bm P}_{A'}^{(2,N)})^{[1,2]}+(1-\gamma)({\bm P}_B^{(2,N)})^{[1,2]},
\end{eqnarray*}
and with $\bm\pi$, $\bm\Pi$, and $\bm Z$ chosen accordingly and $\bm 1:=(1,1,1)^\T$, we can evaluate
\begin{eqnarray*}
\sigma_{(\gamma,1-\gamma)}^{([1,2],N)}&:=&\bm\pi\bm P^{[1,2]}\bm1-(\bm\pi\bm P^{[1]}\bm1)(\bm\pi\bm P^{[2]}\bm1)\\
&&\quad{}+\bm\pi\bm P^{[1]}(\bm Z-\bm\Pi)\bm P^{[2]}\bm 1
+\bm\pi\bm P^{[2]}(\bm Z-\bm\Pi)\bm P^{[1]}\bm 1
\end{eqnarray*}
as a function of $N$, at least if we fix $\rho$ and $\gamma$.

With $\rho=1/3$ and $\gamma=1/2$, we conclude that
\begin{eqnarray}\label{sigmaC2}
(\sigma_{(1/2,1/2)}^{(N)})^2&=&27 (-36821493886409 + 71724260647553 N - 46282959184439 N^2 \nonumber\\
&&\qquad{}+ 9902542819695 N^3)\\
&&\quad/[8331019058 (-269171 + 524347 N - 338381 N^2 + 72405 N^3)],\nonumber
\end{eqnarray}
which is monotonically increasing in $N\ge2$, ranging from
$$
(\sigma_{(1/2,1/2)}^{(2)})^2={114315959583\over258261590798}\approx0.442636
$$
to
$$
\lim_{N\to\infty}(\sigma_{(1/2,1/2)}^{(N)})^2={5941525691817\over13404609664322}\approx0.443245.
$$

Let us summarize our results for random mixtures.
Let $S_n$ be the cumulative profit after $n$ turns to the ensemble of $N\ge2$ players playing the mixed game $\gamma A'+(1-\gamma)B$, where $0\le\gamma\le1$.  We assume the parameterization (\ref{param}) with $\eps=0$.

\begin{theorem}\label{limit-thm}
If $\gamma=1$ so that game $A'$ is always played, then $\P(S_n=0$ for all $n\ge1)=1$.

If $\gamma=0$ so that game $B$ is always played, then $\{S_n-S_{n-1}\}_{n\ge1}$ satisfies the SLLN and the CLT with mean and variance parameters $\mu_B^{(N)}=0$ and $(\sigma_B^{(N)})^2$ as in (\ref{sigmaB2}).

If $0<\gamma<1$ so that both games are played, then $\{S_n-S_{n-1}\}_{n\ge1}$ satisfies the SLLN and the CLT with mean and variance parameters $\mu_{(\gamma,1-\gamma)}^{(N)}$ as in (\ref{mu_C}) (or (\ref{mu_C-2})) and $(\sigma_{(\gamma,1-\gamma)}^{(N)})^2$, at least when $\rho=1/3$ and $\gamma=1/2$, as in (\ref{sigmaC2}).  When $\rho\ne1/3$ or $\gamma\ne1/2$, we implicitly assume that $(\sigma_{(\gamma,1-\gamma)}^{(N)})^2>0$.
\end{theorem}

\begin{proof}
The first conclusion is obvious.  The second and third conclusions follow from Theorem 1, though the mean and variance parameters are obtained not from the theorem but by using the methods described in the text.
\end{proof}

To compare our results with those of Toral (2002), we must restore the bias parameter $\eps>0$.  For simplicity, let us take $\gamma=1/2$, as he did.  Then
\begin{eqnarray}\label{mu_C:eps}
\mu_{(1/2,1/2)}^{(N)}&=&\{3 [2 (1 - \rho)^3 (1 + \rho) - \eps (13 + 26 \rho + 30 \rho^2 + 26 \rho^3 + 13 \rho^4)\\
&&\;{}+ \eps^2 (1 - \rho)^3 (1 + \rho) - 2 \eps^3 (1 + \rho)^2 (1 + \rho^2) ]\}/\{2 [2(10 + 20 \rho \nonumber\\
&&\;{}+ 21 \rho^2 + 20 \rho^3 + 10 \rho^4) - \eps (1 - \rho)^3 (1 + \rho) + 3 \eps^2 (1 + \rho)^2 (1 + \rho^2)]\}.\nonumber
\end{eqnarray}
Toral reported a simulation with $\rho=1/3$, $\gamma=1/2$, $\eps=1/100$, and $N=200$.  Actually, $\eps=1/1000$ was intended (personal communication 2011).  With $\rho=1/3$ and $\eps=1/1000$, (\ref{mu_C:eps}) reduces to $193387599/6704101000\approx0.028846$, with which Toral's estimate, 0.029, is consistent.

\section{Mean profit for nonrandom patterns}\label{patterns}

Toral (2002) omitted discussion of the case in which his games $A'$ and $B$ are played in a nonrandom periodic pattern such as $A'BBA'BBA'BB\cdots$.  Let us denote by $[r,s]$ the pattern $(A')^rB^s$ repeated ad infinitum.  We would like to apply the results of Ethier and Lee (2009) to the pattern $[r,s]$, showing that the Parrondo effect is present for all $r,s\ge1$.  (Unlike in the original one-player Parrondo games, the case $r=s=1$ is included.)  We do this by showing that the mean profit per turn for the ensemble of players, $\mu_{[r,s]}^{(N)}$, is positive if $0<\rho<1$, zero if $\rho=1$, and negative if $\rho>1$, for all $r,s\ge1$ and $N\ge2$.  As we will see, here the mean parameter depends on $N$ and it takes a particularly simple form in the limit as $N\to\infty$.

First, Theorem 6 of Ethier and Lee (2009) is applicable.  (The assumption there that $\bm P_A$ is irreducible and aperiodic is unnecessary.)  But again it is simplest to apply the results to one or two players at a time, as we did in Sections \ref{alternative} and \ref{variance_C}.  Let us begin by finding the mean parameter $\mu_{[r,s]}^{(N)}$.

For the original one-player Parrondo games, in which
$$
\setlength{\arraycolsep}{1mm}
{\bm P}_A:={1\over2}\left(\begin{array}{ccc}
0&1&1\\
1&0&1\\
1&1&0
\end{array}\right),\quad
\setlength{\arraycolsep}{1mm}
{\bm P}_B:=\left(\begin{array}{ccc}
0&p_0&q_0\\
q_1&0&p_1\\
p_2&q_2&0
\end{array}\right),\quad
\setlength{\arraycolsep}{1mm}
{\bm W}:=\left(\begin{array}{rrr}
0&1&-1\\
-1&0&1\\
1&-1&0
\end{array}\right).
$$
Ethier and Lee (2009) showed that
$$
\mu_{[r,s]}={1\over r+s}\,{\bm \pi}_{s,r}{\bm R}\,\,{\rm diag}\!\left(s,\,{1-e_1^s\over1-e_1},\,{1-e_2^s\over1-e_2}\right){\bm L}{\bm\zeta},
$$
where $\bm\pi_{s,r}$ is the unique stationary distribution of $\bm P_B^s\bm P_A^r$, $\bm R$ is the matrix of right eigenvectors of $\bm P_B$, $e_1$ and $e_2$ are the nonunit eigenvalues of $\bm P_B$, $\bm L:=\bm R^{-1}$, and $\bm\zeta:=(\bm P_B\circ\bm W)\bm1$.  They further showed that this formula reduces algebraically to
$$
\mu_{[r,s]}=E_{r,s}/D_{r,s},
$$
where
\begin{eqnarray}\label{Ers}
E_{r,s}&:=&3a_r\{[2+(3a_r-1)(e_1^s+e_2^s-2e_1^s e_2^s)-(e_1^s+e_2^s)](1-\rho)(1+\rho)S\nonumber\\
&&\qquad\qquad{}+a_r(e_2^s-e_1^s)[5(1+\rho)^2(1+\rho^2)-4\rho^2]\}(1-\rho)^2
\end{eqnarray}
and
\begin{equation}\label{Drs}
D_{r,s}:=4(r+s)[1+(3a_r-1)e_1^s][1+(3a_r-1)e_2^s](1+\rho+\rho^2)^2S
\end{equation}
with $a_r:=(1-(-1/2)^r)/3$ and $S:=\sqrt{(1+\rho^2)(1+4\rho+\rho^2)}$.

We apply these results but with $\bm P_A$ and $\bm P_B$ replaced by
$$
\bm P_{A'}^{(1,N)}:=N^{-1}[2\bm P_A^{(1)}+(N-2)\bm I_3]\quad{\rm and}\quad
\bm P_B^{(1,N)}:=N^{-1}[\bm P_B^{(1)}+(N-1)\bm I_3].
$$
Now $(\bm P_{A'}^{(1,N)})^r$ is given by the same formula as $\bm P_A^r$ but with $a_r$ redefined as
\begin{equation}\label{a_r}
a_r:=[1-(1-3/N)^r]/3,
\end{equation}
and $(\bm P_B^{(1,N)})^s$ has the same spectral representation as $\bm P_B^s$ but with the nonunit eigenvalues replaced by
\begin{equation}\label{e1,e2}
e_1:=1-{1-e_1^\circ\over N},\qquad e_2:=1-{1-e_2^\circ\over N},
\end{equation}
where $e_1^\circ$ and $e_2^\circ$ are the nonunit eigenvalues of $\bm P_B$, namely
$$
e_1^\circ:=-{1\over2}+{(1-\rho)S\over2(1+\rho)(1+\rho^2)},\qquad
e_2^\circ:=-{1\over2}-{(1-\rho)S\over2(1+\rho)(1+\rho^2)}.
$$
The matrices $\bm R$ and $\bm L$ are unchanged.

We conclude that
\begin{equation}\label{mu[r,s]-eigenvalues}
\mu_{[r,s]}^{(N)}=NE_{r,s}/D_{r,s},
\end{equation}
where $E_{r,s}$ and $D_{r,s}$ are as in (\ref{Ers}) and (\ref{Drs}) with only the changes (\ref{a_r}) and (\ref{e1,e2}).  For example, this leads to
\begin{eqnarray*}
\mu_{[1,1]}^{(N)}&=&3N (2N-3) (1 - \rho)^3 (1 + \rho)/
 \{2 [18 (1 + \rho + \rho^2)^2 -
    3 N (13 + 26 \rho \\
    &&\quad{}+ 30 \rho^2 + 26 \rho^3 + 13 \rho^4) +
    2N^2 (10 + 20 \rho + 21 \rho^2 + 20 \rho^3 + 10 \rho^4)]\}
\end{eqnarray*}
and
\begin{eqnarray*}
&&\!\!\!\!\!\mu_{[1,2]}^{(N)}\\
&&{}=2 N (1 - \rho)^3 (1 + \rho) [- 3 (1 + \rho + \rho^2)^2 + N (10 + 20 \rho + 21 \rho^2 + 20 \rho^3 + 10 \rho^4)\\
&&\quad{} -9 N^2 (1 + \rho)^2 (1 + \rho^2) + 3 N^3 (1 + \rho)^2 (1 +\rho^2)]
/[36 (1 + \rho + \rho^2)^4\\
&&\quad{} - 12 N (1 + \rho + \rho^2)^2 (11 + 22 \rho + 24 \rho^2 + 22 \rho^3 + 11 \rho^4) + N^2 (193 + 772 \rho\\
&&\quad{} + 1660 \rho^2 + 2548 \rho^3 + 2938 \rho^4 + 2548 \rho^5 + 1660 \rho^6 + 772 \rho^7 + 193 \rho^8) \\
&&\quad{} - 3 N^3 (1 + \rho)^2 (43 + 86 \rho + 145 \rho^2 + 172 \rho^3 + 145 \rho^4 + 86 \rho^5 + 43 \rho^6) \\
&&\quad{} + N^4 (1 + \rho)^2 (35 + 70 \rho + 113 \rho^2 + 140 \rho^3 + 113 \rho^4 + 70 \rho^5 + 35 \rho^6)].
\end{eqnarray*}
Both of these functions are positive for all $N\ge2$, as can be seen by expanding numerators and denominators in powers of $N-2$ and noticing that all coefficients are polynomials in $\rho$ with only positive coefficients.

Although these formulas become increasingly complicated as $r$ and $s$ increase, their limits as $N\to\infty$ have a very simple form.  To see this, it suffices to note that
$$
a_r={r\over N}+O\bigg({1\over N^2}\bigg),\qquad e_1^s=1-{(1-e_1^\circ)s\over N}+O\bigg({1\over N^2}\bigg),
$$
and similarly for $e_2^s$, so (\ref{mu[r,s]-eigenvalues}) converges as $N\to\infty$ to
$$
{3 r s (1 - \rho)^3 (1 + \rho) \over9 r^2 (1 + \rho)^2 (1 + \rho^2) + 9 r s (1 + \rho)^2 (1 + \rho^2) + 2 s^2 (1 + \rho + \rho^2)^2},
$$
which coincides with (\ref{mu_C}) (or (\ref{mu_C-2})) when $\gamma=r/(r+s)$.  This limit is positive if $0<\rho<1$, zero if $\rho=1$, and negative if $\rho>1$, so we conclude that the Parrondo effect is present for all $r,s\ge1$, as long as $N$ is large enough and $\rho\ne1$.  This relationship between the random-mixture case and the nonrandom-pattern case is not present in the original one-player Parrondo games except in a single case ($r=2$, $s=1$).  (We have confirmed this for $r,s\ge1$ and $r+s\le75$ and expect that it is true generally.)

We now verify that the Parrondo effect is always present. We begin with a lemma.

\begin{lemma}\label{lemma1}
If $0<a<b<c$, then $(c^n-b^n)/(b^n-a^n)$ is increasing in $n\ge1$.
\end{lemma}

\begin{proof}
Divide both numerator and denominator by $b^n$ to see that we can, without loss of generality, assume that $b=1$.  So the aim is to show that
$$
{c^n-1\over1-a^n}<{c^{n+1}-1\over1-a^{n+1}}, \qquad n\ge1,
$$
or that
$$
{c^n-1\over c^{n+1}-1}<{a^n-1\over a^{n+1}-1}, \qquad n\ge1.
$$
For this it is enough to fix $n\ge1$ and show that the function
$$
f(x):={x^n-1\over x^{n+1}-1},
$$
defined by continuity at $x=1$, is decreasing on $(0,\infty)$.  Its derivative has the same sign as
$$
-[x^{n+1}-(n+1)x+n],
$$
so it is enough that the quantity within brackets is positive for $x>1$ and $0<x<1$.  First suppose that $x>1$.  Then
\begin{eqnarray*}
x^{n+1}-(n+1)x+n&=&(x-1+1)^{n+1}-(n+1)(x-1)-1\\
&=&(x-1)^{n+1}+{n+1\choose1}(x-1)^n+\cdots+{n+1\choose n-1}(x-1)^2\\
&>&0.
\end{eqnarray*}
Next suppose that $0<x<1$.  Then
\begin{eqnarray*}
x^{n+1}-(n+1)x+n&=&x^{n+1}-1-(n+1)(x-1)\\
&=&(x-1)(x^n+x^{n-1}+\cdots+x+1)-(n+1)(x-1)\\
&=&(x-1)[x^n+x^{n-1}+\cdots+x+1-(n+1)]\\
&>&0.
\end{eqnarray*}
This completes the proof.
\end{proof}

\begin{theorem}\label{positivity}
$\mu_{[r,s]}^{(N)}$ is positive if $0<\rho<1$, zero if $\rho=1$, and negative if $\rho>1$, for all $r,s\ge1$ and $N\ge2$.
\end{theorem}

\begin{proof}
Denoting $\mu_{[r,s]}^{(N)}$ temporarily by $\mu_{[r,s]}^{(N)}(\rho)$ to emphasize its dependence on $\rho$, it can be shown algebraically or probabilistically that
$$
\mu_{[r,s]}^{(N)}(1/\rho)=-\mu_{[r,s]}^{(N)}(\rho),
$$
so it will suffice to treat the case $0<\rho<1$.
First, $|3a_r-1|<1$ and $e_1,e_2\in(0,1)$, so $D_{r,s}>0$.  Since $a_r>0$, it suffices to show that
\begin{eqnarray*}
&&[2+(3a_r-1)(e_1^s+e_2^s-2e_1^s e_2^s)-(e_1^s+e_2^s)](1-\rho)(1+\rho)S\\
&&\qquad\qquad{}+a_r(e_2^s-e_1^s)[5(1+\rho)^2(1+\rho^2)-4\rho^2]>0.
\end{eqnarray*}
Discarding the $-4\rho^2$ term (since $e_2^s-e_1^s<0$), it is enough to show that
\begin{eqnarray}\label{suff}
&&(1-e_1^s)[1+(3a_r-1)e_2^s]+(1-e_2^s)[1+(3a_r-1)e_1^s]\nonumber\\
&&\qquad\qquad{}-a_r(e_1^s-e_2^s){5(1+\rho)(1+\rho^2)\over(1-\rho)S}>0.
\end{eqnarray}
Now $e_1^\circ=(-1+x)/2$ and $e_2^\circ=(-1-x)/2$, where $x:=(1-\rho)S/[(1+\rho)(1+\rho^2)]\in(0,1)$, so
$e_1=(2N-3+x)/(2N)$ and $e_2=(2N-3-x)/(2N)$.

Let us first assume that $N\ge3$.  Then $3a_r-1\le0$, so, replacing $e_1^s$ and $e_2^s$ within brackets in (\ref{suff}) by 1, we need only show that
$$
3(1-e_1^s)+3(1-e_2^s)>(e_1^s-e_2^s){5(1+\rho)(1+\rho^2)\over(1-\rho)S},
$$
or that
\begin{equation}\label{x-ineq}
{2(2N)^s-[(2N-3+x)^s+(2N-3-x)^s]\over[(2N-3+x)^s-(2N-3-x)^s]/x}>{5\over3}.
\end{equation}
The denominator is a polynomial of degree $s-1$ in $x$ with positive coefficients, while the term within brackets in the numerator is a polynomial of degree $s$ in $x$ with positive coefficients.  So the left side of (\ref{x-ineq}) is decreasing in $x$, and it suffices to verify it at $x=1$.  For this we notice that
\begin{eqnarray*}
{2(2N)^s-[(2N-2)^s+(2N-4)^s]\over(2N-2)^s-(2N-4)^s}=2{N^s-(N-1)^s\over(N-1)^s-(N-2)^s}+1,
\end{eqnarray*}
and the fraction on the right is increasing in $s\ge1$ by Lemma \ref{lemma1}.  At $s=1$ the value is 3, so the desired inequality holds.

It remains only to consider the case $N=2$.  The same argument works if $r$ is even because then $3a_r-1\le0$ still holds.  If $r$ is odd, we can replace the quantities within brackets in (\ref{suff}) by 1 and can replace $a_r$ in the second line of (\ref{suff}) by $a_1=1/2$.  Thus, we need only verify (\ref{x-ineq}) with 5/3 replaced by 5/2, and of course it still holds.
\end{proof}

\section{Remark on a ``coincidence''}\label{coincidence}

We can prove algebraically that
\begin{eqnarray*}
\lim_{M\to\infty}\mu_{[r,r]}^{(M)}&=&\mu_{(1/2,1/2)}^{(N)}=(3/2)\mu_{(2/3,1/3)}^{(1)}=(3/2)\mu_{[2,1]}^{(1)}=\mu_{[1,1]}^{(2)}\\
&=&{3 (1 - \rho)^3 (1 + \rho)\over 2(10 + 20 \rho + 21 \rho^2 + 20 \rho^3 + 10 \rho^4)}
\end{eqnarray*}
for all $r\ge1$ and $N\ge2$, where superscripts refer to the number of players.  (For superscripts equal to 1, the games are $A$ and $B$, the original one-player Parrondo games.  For superscripts 2 or larger, the games are $A'$ and $B$.)  The first equality is from Section \ref{patterns}.  Can the others be explained probabilistically?  We can elucidate at least the second equality, while the third and fourth remain partially unexplained.

Since $\mu_{(\gamma,1-\gamma)}^{(N)}$ does not depend on $N$, it is enough to verify the identity with $N=2$.  Let us consider the profit of one player when two players are playing.  Recalling (\ref{pB1N}) and (\ref{pA1N}) with $N=2$, we have
$$
\bm P_{A'}^{(1,2)}:=\bm P_A^{(1)}\qquad{\rm and}\qquad\bm P_B^{(1,2)}:=(1/2)(\bm P_B^{(1)}+\bm I_3).
$$
The former is just the one-step transition matrix for the original one-player game $A$, and we have
$$
(1/2)\bm P_{A'}^{(1,2)}+(1/2)\bm P_B^{(1,2)}=(1/2)\bm P_A^{(1)}+(1/4)\bm P_B^{(1)}+(1/4)\bm I_3.
$$
The left side describes the $(1/2,1/2)$ mixture of games $A'$ and $B$, as viewed by one of two players.  Its mean is $(1/2)\mu_{(1/2,1/2)}^{(2)}$.  The right side describes the $(2/3,1/3)$-mixture of games $A$ and $B$ if we ignore the $(1/4)\bm I_3$ term and normalize to ensure a stochastic matrix.  That term just slows down the process, making one-fourth of its transitions null.  So its mean is $(3/4)\mu_{(2/3,1/3)}^{(1)}$.  These are equal, so $\mu_{(1/2,1/2)}^{(2)}=(3/2)\mu_{(2/3,1/3)}^{(1)}$, as claimed.

This can be regarded as a more correct version of the argument sketched in the third paragraph of page L307 of Toral (2002) and attributed to an anonymous referee of that paper.

\section{Variance parameter for nonrandom patterns}

We can evaluate the variance parameter $(\sigma_{[r,s]}^{(N)})^2$ for the nonrandom pattern $[r,s]$ in the $N$-player games directly for small $N$, using the state space $\bar\Sigma_N$ with its ${N+2\choose2}$ states.  We apply (25)--(27) of Ethier and Lee (2009), obtaining, for example,
\begin{eqnarray*}
&&(\sigma_{[1,1]}^{(2)})^2\\
&&\quad{}=[9 (466 + 2680 \rho + 7621 \rho^2 + 16310 \rho^3 + 29018 \rho^4 + 41582 \rho^5 + 51471 \rho^6\\
&&\qquad\;{} + 55998 \rho^7 + 51471 \rho^8 + 41582 \rho^9 + 29018 \rho^{10} + 16310 \rho^{11} + 7621 \rho^{12}\\
&&\qquad\;{} + 2680 \rho^{13} + 466 \rho^{14})]/[4 (2 - \rho + 2 \rho^2) (10 + 20 \rho + 21 \rho^2 + 20 \rho^3 + 10 \rho^4)^3],
\end{eqnarray*}
which reduces when $\rho=1/3$ to $74176355601/141627323986\approx0.523743$.  Since $N=2$, this is a computation involving $6\times6$ matrices.

To get results for larger $N$, we apply the method of considering one or two players at a time.  By analogy with (\ref{sigma_C^N}), we have
\begin{equation*}
(\sigma_{[r,s]}^{(N)})^2=N(\sigma_{[r,s]}^{(1,N)})^2+N(N-1)\sigma_{[r,s]}^{([1,2],N)},
\end{equation*}
where
\begin{eqnarray*}
(\sigma_{[r,s]}^{(1,N)})^2&=&{1\over r+s}\bigg\{\sum_{u=0}^{r-1}[\bm\pi\bm P_A^u\ddot{\bm P}_A\bm1-(\bm\pi\bm P_A^u\dot{\bm P}_A\bm1)^2]\\
&&\qquad\quad{} +\sum_{v=0}^{s-1}[\bm\pi\bm P_A^r\bm P_B^v\ddot{\bm P}_B\bm1-(\bm\pi\bm P_A^r\bm P_B^v\dot{\bm P}_B\bm1)^2] \nonumber\\
&&\qquad\quad{}+2\sum_{0\le u<v\le r-1}\bm\pi\bm P_A^u\dot{\bm P}_A(\bm P_A^{v-u-1}-\bm\Pi\bm P_A^v)\dot{\bm P}_A\bm1 \nonumber\\
&&\qquad\quad{}+2\sum_{u=0}^{r-1}\sum_{v=0}^{s-1}\bm\pi\bm P_A^u\dot{\bm P}_A (\bm P_A^{r-u-1}-\bm\Pi\bm P_A^r)\bm P_B^v\dot{\bm P}_B\bm1\nonumber\\
&&\qquad\quad{}+2\sum_{0\le u<v\le s-1}\bm\pi\bm P_A^r\bm P_B^u\dot{\bm P}_B(\bm P_B^{v-u-1}-\bm\Pi\bm P_A^r\bm P_B^v)\dot{\bm P}_B\bm1\\
&&\qquad\quad{}+2\bigg[\sum_{u=0}^{r-1}\sum_{v=0}^{r-1}\bm\pi\bm P_A^u\dot{\bm P}_A\bm P_A^{r-u-1}\bm P_B^s(\bm Z-\bm\Pi)\bm P_A^v\dot{\bm P}_A\bm1 \nonumber\\
&&\qquad\qquad\quad{}+\sum_{u=0}^{r-1}\sum_{v=0}^{s-1}\bm\pi\bm P_A^u\dot{\bm P}_A\bm P_A^{r-u-1}\bm P_B^s(\bm Z-\bm\Pi)\bm P_A^r\bm P_B^v\dot{\bm P}_B\bm1 \nonumber\\
&&\qquad\qquad\quad{}+\sum_{u=0}^{s-1}\sum_{v=0}^{r-1}\bm\pi\bm P_A^r\bm P_B^u\dot{\bm P}_B\bm P_B^{s-u-1}(\bm Z-\bm\Pi)\bm P_A^v\dot{\bm P}_A\bm1 \nonumber\\
&&\qquad\qquad\quad{}+\sum_{u=0}^{s-1}\sum_{v=0}^{s-1}\bm\pi\bm P_A^r\bm P_B^u\dot{\bm P}_B\bm P_B^{s-u-1}(\bm Z-\bm\Pi)\bm P_A^r\bm P_B^v\dot{\bm P}_B\bm1\bigg]\bigg\}
\end{eqnarray*}
(from Ethier and Lee 2009) with $\bm P_A$ and $\bm P_B$ replaced by $\bm P_{A'}^{(1,N)}$ and $\bm P_B^{(1,N)}$ as defined in Section \ref{alternative}.

The covariance term, $\sigma_{[r,s]}^{([1,2],N)}$, requires an extension of the preceding formula to covariances.  We omit the details of the derivation and just give the result:
\begin{eqnarray*}
\sigma_{[r,s]}^{([1,2],N)}&=&{1\over r+s}\bigg\{\sum_{u=0}^{r-1}[\bm\pi\bm P_A^u\bm P_A^{[1,2]}\bm1-(\bm\pi\bm P_A^u\bm P_A^{[1]}\bm1)(\bm\pi\bm P_A^u\bm P_A^{[2]}\bm1)]\\
&&\quad\qquad{}+\sum_{v=0}^{s-1}[\bm\pi\bm P_A^r\bm P_B^v\bm P_B^{[1,2]}\bm1-(\bm\pi\bm P_A^r\bm P_B^v\bm P_B^{[1]}\bm1)(\bm\pi\bm P_A^r\bm P_B^v\bm P_B^{[2]}\bm1)]\\
&&\quad\qquad{}+\sum_{0\le u<v\le r-1}[\bm\pi\bm P_A^u\bm P_A^{[1]}(\bm P_A^{v-u-1}-\bm\Pi\bm P_A^v)\bm P_A^{[2]}\bm1\\
\noalign{\vglue-3mm}
&&\qquad\qquad\qquad\qquad\qquad{}+\bm\pi\bm P_A^u\bm P_A^{[2]}(\bm P_A^{v-u-1}-\bm\Pi\bm P_A^v)\bm P_A^{[1]}\bm1]\\
&&\quad\qquad{}+\sum_{u=0}^{r-1}\sum_{v=0}^{s-1}[\bm\pi\bm P_A^u\bm P_A^{[1]}(\bm P_A^{r-u-1}-\bm\Pi\bm P_A^r)\bm P_B^v\bm P_B^{[2]}\bm1\\
\noalign{\vglue-2mm}
&&\qquad\qquad\qquad\qquad\qquad{}+\bm\pi\bm P_A^u\bm P_A^{[2]}(\bm P_A^{r-u-1}-\bm\Pi\bm P_A^r)\bm P_B^v\bm P_B^{[1]}\bm1]\\
&&\quad\qquad{}+\sum_{0\le u<v\le s-1}[\bm\pi\bm P_A^r\bm P_B^u\bm P_B^{[1]}(\bm P_B^{v-u-1}-\bm\Pi\bm P_A^r\bm P_B^v)\bm P_B^{[2]}\bm1\\
\noalign{\vglue-3mm}
&&\qquad\qquad\qquad\qquad\qquad{}+\bm\pi\bm P_A^r\bm P_B^u\bm P_B^{[2]}(\bm P_B^{v-u-1}-\bm\Pi\bm P_A^r\bm P_B^v)\bm P_B^{[1]}\bm1]\\
&&\quad\qquad{}+\sum_{u=0}^{r-1}\sum_{v=0}^{r-1}[\bm\pi\bm P_A^u\bm P_A^{[1]}\bm P_A^{r-u-1}\bm P_B^s(\bm Z-\bm\Pi)\bm P_A^v\bm P_A^{[2]}\bm1\\
\noalign{\vglue-3mm}
&&\qquad\qquad\qquad\qquad{}+\bm\pi\bm P_A^u\bm P_A^{[2]}\bm P_A^{r-u-1}\bm P_B^s(\bm Z-\bm\Pi)\bm P_A^v\bm P_A^{[1]}\bm1]\\
&&\quad\qquad{}+\sum_{u=0}^{r-1}\sum_{v=0}^{s-1}[\bm\pi\bm P_A^u\bm P_A^{[1]}\bm P_A^{r-u-1}\bm P_B^s(\bm Z-\bm\Pi)\bm P_A^r\bm P_B^v\bm P_B^{[2]}\bm1\\
\noalign{\vglue-3mm}
&&\qquad\qquad\qquad\qquad{}+\bm\pi\bm P_A^u\bm P_A^{[2]}\bm P_A^{r-u-1}\bm P_B^s(\bm Z-\bm\Pi)\bm P_A^r\bm P_B^v\bm P_B^{[1]}\bm1]\\
&&\quad\qquad{}+\sum_{u=0}^{s-1}\sum_{v=0}^{r-1}[\bm\pi\bm P_A^r\bm P_B^u\bm P_B^{[1]}\bm P_B^{s-u-1}(\bm Z-\bm\Pi)\bm P_A^v\bm P_A^{[2]}\bm1\\
\noalign{\vglue-3mm}
&&\qquad\qquad\qquad\qquad{}+\bm\pi\bm P_A^r\bm P_B^u\bm P_B^{[2]}\bm P_B^{s-u-1}(\bm Z-\bm\Pi)\bm P_A^v\bm P_A^{[1]}\bm1]\\
&&\quad\qquad{}+\sum_{u=0}^{s-1}\sum_{v=0}^{s-1}[\bm\pi\bm P_A^r\bm P_B^u\bm P_B^{[1]}\bm P_B^{s-u-1}(\bm Z-\bm\Pi)\bm P_A^r\bm P_B^v\bm P_B^{[2]}\bm1\\
\noalign{\vglue-3mm}
&&\qquad\qquad\qquad\qquad{}+\bm\pi\bm P_A^r\bm P_B^u\bm P_B^{[2]}\bm P_B^{s-u-1}(\bm Z-\bm\Pi)\bm P_A^r\bm P_B^v\bm P_B^{[1]}\bm1]\bigg\}
\end{eqnarray*}
with $\bm P_A$ and $\bm P_B$ replaced by $\bm P_{A'}^{(2,N)}$ and $\bm P_B^{(2,N)}$ as defined in Section \ref{variance_C}.

By analogy to (\ref{sigmaC2}), with $\rho=1/3$, we conclude that
\begin{eqnarray}\label{sigma[1,1]2}
(\sigma_{[1,1]}^{(N)})^2&=&9 (615639408424560 - 6408926620214040 N + 29541545957894139 N^2\nonumber \\
&&\quad{} - 80214814200037491 N^3 + 143582273075781927 N^4\nonumber\\
&&\quad{} - 179192557802543130 N^5 + 160434481099881996 N^6 \nonumber\\
&&\quad{}- 104152159483211664 N^7 + 48799091685478468 N^8\nonumber\\
&&\quad{} - 16137521956595246 N^9 + 3584898779152593 N^{10} \nonumber\\
&&\quad{} - 481633399018397 N^{11} + 29679648590925 N^{12})\nonumber\\
&&\;\;/[2 (1521 - 3174 N + 1609 N^2)^3 (3285360 -
     9816120 N + 12525387 N^2\nonumber\\
&&\quad{} - 8725589 N^3 + 3501928 N^4 - 768851 N^5 +
     72405 N^6)],
\end{eqnarray}
which is monotonically decreasing in $N$, ranging from
$$
(\sigma_{[1,1]}^{(2)})^2={74176355601\over141627323986}\approx0.523743
$$
to
$$
\lim_{N\to\infty}(\sigma_{[1,1]}^{(N)})^2={5935929718185\over13404609664322}\approx0.442827.
$$
This last number differs slightly from the corresponding limit in the random-mixture case, showing that the variances lack the nice property that the means enjoy.

Similar formulas can be found for other $[r,s]$, assuming $\rho=1/3$.  With $[r,s]=[1,2]$ (resp., $[2,1]$, $[2,2]$, $[2,4]$, $[3,3]$, $[4,2]$), we get in place of (\ref{sigma[1,1]2}) a ratio of degree-21 (resp., 24, 33, 51, 54, 57) polynomials in $N$, and
$$
\lim_{N\to\infty}(\sigma_{[r,s]}^{(N)})^2=\begin{cases}
{1891312136577\over6060711605323}\approx0.312061&\text{if $[r,s]=[2,1]$ or $[4,2]$,}\\
\noalign{\medskip}
{5935929718185\over13404609664322}\approx0.442827&\text{if $[r,s]=[1,1]$ or $[2,2]$ or $[3,3]$,}\\
\noalign{\medskip}
{136286243910\over252688187761}\approx0.539346&\text{if $[r,s]=[1,2]$ or $[2,4]$.}
\end{cases}
$$
In particular, it appears that the result for $[r,s]$ depends on $r$ and $s$ only through $r/(r+s)$.  We have confirmed this only in the several cases shown above; a proof for all integers $r,s\ge1$ seems difficult.

Let us summarize our results for nonrandom patterns.
Let $S_n$ be the cumulative profit after $n$ turns to the ensemble of $N\ge2$ players playing the nonrandom pattern $(A')^r B^s$ (denoted by $[r,s]$) repeatedly, with $r$ and $s$ being positive integers.  We assume the parameterization (\ref{param}) with $\eps=0$.

\begin{theorem}
$\{S_n-S_{n-1}\}_{n\ge1}$ satisfies the SLLN and the CLT with mean and variance parameters
$\mu_{[r,s]}^{(N)}$ computable for specified $r,s\ge1$ as a function of $N\ge2$ and $\rho>0$ and $(\sigma_{[r,s]}^{(N)})^2$ computable for specified $r,s\ge1$, $N\ge2$, and $\rho>0$.  We implicitly assume that $(\sigma_{[r,s]}^{(N)})^2>0$.
\end{theorem}

\begin{proof}
The proof is as for Theorem \ref{limit-thm}, except that, instead of citing Theorem 1, we cite Theorem 6 of Ethier and Lee (2009).
\end{proof}

\section{Sample variance of players' capitals}

Recall our notation in which $S_n^{[i]}$ is the capital of player $i$ (one of the $N$ players) after $n$ trials, so that
$$
S_n:=\sum_{i=1}^N S_n^{[i]}
$$
is the total capital of the $N$ players after $n$ trials.
Toral (2002) simulated the sample variance of $S_n^{[1]},\ldots,S_n^{[N]}$, namely
$$
{1\over N}\bigg[\sum_{i=1}^N (S_n^{[i]})^2-N\bigg({1\over N}\sum_{i=1}^N S_n^{[i]}\bigg)^2\bigg],
$$
finding that it grows approximately linearly in $n$.  Let us replace this sample variance by its unbiased (at least in the case of a random sample) version,
$$
(\hat{\sigma}^{(N)})_n^2:={1\over N-1}\bigg[\sum_{i=1}^N (S_n^{[i]})^2-N\bigg({1\over N}\sum_{i=1}^N S_n^{[i]}\bigg)^2\bigg],
$$
and let us consider its expectation $\E[(\hat{\sigma}^{(N)})_n^2]$ instead of the random variable itself.  We can evaluate this using the results already obtained.  Indeed,
\begin{eqnarray*}
\E[(\hat{\sigma}^{(N)})_n^2]&:=&{1\over N-1}\bigg[\sum_{i=1}^N \E[(S_n^{[i]})^2]-{1\over N}\E[(S_n)^2]\bigg]\\
&\;=&{1\over N-1}\bigg[\sum_{i=1}^N \{\Var(S_n^{[i]})+(\E[S_n^{[i]}])^2\}-{1\over N}\{\Var(S_n)+(\E[S_n])^2\}\bigg]\\
&\;=&{1\over N-1}\bigg[N\Var(S_n^{[1]})-{1\over N}\Var(S_n)\bigg]\\
&\;\sim&n{1\over N-1}\bigg[N(\sigma^{(1,N)})^2-{1\over N}[N(\sigma^{(1,N)})^2+N(N-1)\sigma^{([1,2],N)}]\bigg]\\
&\;=&n[(\sigma^{(1,N)})^2-\sigma^{([1,2],N)}]
\end{eqnarray*}
as $n\to\infty$, hence
\begin{equation*}
\lim_{n\to\infty}n^{-1}\E[(\hat{\sigma}^{(N)})_n^2]=(\sigma^{(1,N)})^2-\sigma^{([1,2],N)}.
\end{equation*}
We have omitted subscripts $A'$, $B$, $(\gamma,1-\gamma)$, and $[r,s]$ because we want to apply this formula in all cases.

Let us first consider the random-mixture case with $\rho$ arbitrary.  With $\gamma=1/2$ we have
\begin{eqnarray}\label{asymp-sample-var}
&&\!\!\!\!\!\!\!\!\!\!(\sigma_{(1/2,1/2)}^{(1,N)})^2-\sigma_{(1/2,1/2)}^{([1,2],N)}\nonumber\\
&&{}\sim27 (97 + 606 \rho + 1926 \rho^2 + 4262 \rho^3 + 7284 \rho^4 + 9894 \rho^5 + 10911 \rho^6\nonumber \\
&&\qquad\quad{}+ 9894 \rho^7 + 7284 \rho^8 + 4262 \rho^9 + 1926 \rho^{10} + 606 \rho^{11} + 97 \rho^{12})\nonumber\\
&&\qquad{}/[2 N (10 + 20 \rho + 21 \rho^2 + 20 \rho^3 + 10 \rho^4)^3]
\end{eqnarray}
as $N\to\infty$, which is (\ref{sigma_C-asymp}) with $\gamma=1/2$.  With $\gamma=0$ (only game $B$ is played) we have
$$
(\sigma_B^{(1,N)})^2-\sigma_B^{([1,2],N)}=\bigg({3\rho\over1+\rho+\rho^2}\bigg)^2{1\over N}.
$$
Finally, with $\gamma=1$ (only game $A'$ is played) we have
$$
(\sigma_{A'}^{(1,N)})^2-\sigma_{A'}^{([1,2],N)}={2\over N-1}\sim{2\over N}.
$$
As Toral found, the result for the mixed game lies between those for games $A'$ and $B$, and this is true regardless of $\rho>0$.  (Our results are not directly comparable to his because we have taken the bias parameter $\eps$ to be 0.)

Finally, let us consider the nonrandom pattern $[1,1]$.  We find that
$(\sigma_{[1,1]}^{(1,N)})^2-\sigma_{[1,1]}^{([1,2],N)}$ has the same asymptotic value as in (\ref{asymp-sample-var}),
so we have another coincidence.  It appears that these expected sample variances have essentially the same property that the means enjoy, namely that their asymptotic value for the case of the nonrandom pattern $[r,s]$ depends only on $r/(r+s)$ and $\rho$ and is equal to the asymptotic value for the case of the random mixture with $\gamma=r/(r+s)$.  We have confirmed this in the six cases $r,s\ge1$ and $r+s\le4$, but a proof for all integers $r,s\ge1$ seems difficult.

\end{document}